\documentstyle[12pt]{article}

\newcommand{\beq}{\begin{equation}}
\newcommand{\eeq}{\end{equation}}
\newcommand{\beeq}{\begin{eqnarray}}
\newcommand{\eeeq}{\end{eqnarray}}
\newcommand{\bdm}{\begin{displaymath}}
\newcommand{\edm}{\end{displaymath}}

\begin{document}
\thispagestyle{empty}
\begin{flushright}
\vspace*{3cm}
{\large SPBU-IP-94-22}
\end{flushright}
\vspace*{2cm}

\begin{center}
{\large {\bf q-Bosonization of the quantum group} \\
{\bf $ GL_q(2) $ based on the Gauss decomposition}} \\[1cm]

{\bf E.V.Damaskinsky} \footnote
 {  Supported by Russian Foundation for
 Fundamental Research, Grant N 94-01-01157-a.}
 \footnote{ E-mail address:
 zhel@vici.spb.su } \\
 Institute of Military Constructing Engineering \\
 Zacharievskaya st 22, \\
 191194, St Petersburg, Russia \\[1cm]

{\bf M.A.Sokolov} \footnote{ E-mail adress:
 sokol@pmash.spb.su} \\
 St Petersburg Institute of Mashine Building \\
 Poliustrovskii pr 14, \\
 195108, St Petersburg, Russia \\[2cm]

  {\bf Abstract.}
\end{center}
The new method of q-bosonization for quantum groups based on the Gauss
decomposition of a transfer matrix of generators is suggested. The
simplest example of the quantum group $GL_q(2)$ is considered in some
details.

\newpage

{\bf 1}. The method of construction of Lie group and algebra representations
using the realization  their generators in terms of boson oscillator
creation and annihilation operators  is well known. The example of such
construction gives the famous Schwinger representation (bosonization) for
angular momentum operators. For investigations of representations of quantum
groups and algebras [1,2] the similar construction is also useful. However
in this case it is natural to replace the usual boson oscillators on
deformed ones [3-6]. This procedure is called $q$-{\em bosonization. }

For the simple quantum algebras of all classical series from Cartan list  $q$
-bosonization procedure is described with sufficient completeness in [7]
(see also [6] for the noncompact $su_q(1,1)$ and [8-10] for super case).
However, the situation for quantum groups is quite different. It seems first
attempts to this end have been done for the quantum group $GL_q(n)$ in the
special case $q^n=1$ in [11,12]. In the general case, $q\in C\backslash \{0\}
$, the examples of q-bosonization were given in [13] and [14] for the $
GL_q(2)$ and $GL_q(3)$ respectively. The attempt to generalize these results
to the $GL_q(n)$ was undertaken in somewhat complicated fashion in Ref.
[15]. The receipt given in [15] based on the special construction of the
representation in the q-analog of the highest weight module. The application
of this receipt to the $n>3$ case needs very tedious computations. The
deficiency of this approach [13,15] consists in the use of specific features
of the Fock representation for the q-oscillators. As a result, a wide class
of non Fock representations (listed, for instance, in [5,16-17]) are
excluded from considerations {\it ab initio}. We remark also that the direct
ways to extend this receipt to quantum groups of other series are absent.
Let us also mention the interesting works [18,19] concerned with the similar
problem for the matrix pseudogroups $S_\mu U(n)$ and used somewhat different
form of q-oscillators. Thus we see that the search of new ways of
realization of $q$-bosonization procedure seems very desirable.

In this paper we want to show that the Gauss decomposition [20] for the $
GL_q(2)$ suggests q-bosonization in a very simple and pure algebraic way.
Let us note that similar methods applied to the general case of the $GL_q(n)$
quantum group [21] as well as to quantum groups of $B_n,C_n$ and $D_n$
series [22]

The Russian version of this work submited to Theor.Math.Phys. in April,
1995.

{\bf 2}. Let us recall [1] that quantum group $GL_q(2)$ is defined as
associative unital algebra freely generated by four generators $a,b,c,d$
(usually written as entries of the matrix $T=\left(
\begin{array}{cc}
a & b \\
c & d
\end{array}
\right) $) factorized by 2-ideal generated by commutation relations of the
form
\begin{eqnarray}\begin{array}{c}
ab=qba,\quad ac=qca,\quad bd=qdb,\quad cd=qdc \\
bc=cb,\quad ad-da=\lambda bc,
\end{array}
\end{eqnarray}where $\lambda \equiv q-q^{-1}$. From this commutation
relations it follows that the quantum ($q$-)determinant
$$
{\cal D}\equiv {\cal D}_q(T)\equiv \det \nolimits_q(t):=ad-qbc=da-q^{-1}bc
$$
commutes with all generators and thus belongs to the center of quantum group
$GL_q(2)$. Let us suppose that ${\cal D}_q(T)\equiv \det \nolimits_q(t)\neq
0 $. The additional condition ${\cal D}_q(T)=1$ pick out the quantum group $
SL_q(2)$. Let us recall that the quantum groups $GL_q(2)$ and its subgroup $
SL_q(2)$ are not groups in usual sense. In particular, the entries of the
matrix $T^n$ subject to the modificated commutation relations (1) with $q$
replaced by $q^{n}$.

{\bf 3}. {\em Quantum deformed oscillator }( or simply $q$-{\em oscillator})
is [3-5] an associative algebra {\cal A$_q$} generated by three elements $
a(=a_{-})$, $a^{\dagger }(=a_{+})$ and $N$, which fulfill the following
commutation relations
\begin{equation}
[a,a^{\dagger }]_q\equiv aa^{\dagger }-qa^{\dagger }a=q^{-N},\quad
[N,a]=-a,\quad [N,a^{\dagger }]=a^{\dagger }.
\end{equation}
{}From this relations it follows that $[N,a^{\dagger }a]=0=[N,aa^{\dagger }]$.
The element $\zeta =q^{-N}([N]-$ $a^{\dagger }a)$ , where $[N]=\frac{
q^N-q^{-N}}{q-q^{-1}}$, commutes with every generator and thus belong to the
center of the $q$-oscillator algebra {\cal A$_q$}. It is worth-while to
stress that the deformed boson oscillator algebra has nontrivial center
unlike the undeformed case. This nontriviality of the center is one of the
reasons for existence of the non-Fock representations of the $q$-oscillator
algebra {\cal A$_q$} [5,16,17]. The $q$-analog of the Fock representation
can be constructed easily from the usual Fock representation of boson
oscillator. Let \{$|n>:n=0,1,2,...\}$ be the standard basis in the Fock
space {\cal H} with vacuum $|0>$. Then the operators $a,a^{\dagger }$ and $N$
, defined by the following action on this basis
$$
N|n>=n|n>,\quad a|n>=\sqrt{[n]}|n-1>,\quad a^{\dagger }|n>=\sqrt{[n+1]}
|n+1>,
$$
fulfill the commutation relations (3) and define the Fock representation for
$q$-oscillator algebra {\cal A$_q$} . As in the usual case the basic vector $
|n>$ can be obtained from the vacuum state by repeated action of the $q$
-creation operator $a^{\dagger }$, namely, $|n>=([n]!)^{-1/2}(a^{\dagger
})^n|0>$. On the Fock space $q$-oscillator operators $a,a^{\dagger }$ are
connected with usual boson oscillator operators $b,b^{\dagger }$ by the
relations
$$
a=b([N]/N)^{1/2},\quad a^{\dagger }=([N]/N)^{1/2}b^{\dagger },\quad
N=b^{\dagger }b.
$$
{}From these relations it follows that Fock basis is not deformed (but the $q$
-analogs of coherent states was subjects to deformation [5]). Now, let us
note that on the Fock space {\cal H} the central element $\zeta $ vanishes ($
\zeta =0)$ . From this fact it follows that on the Fock space {\cal H}
special relations
\begin{equation}
[N]=a^{\dagger }a,\quad [N+1]=aa^{\dagger }
\end{equation}
are valid, as well as the additional commutation relation
\begin{equation}
aa^{\dagger }-q^{-1}a^{\dagger }a=q^N,
\end{equation}
which realize the specific $q\leftrightarrow q^{-1}$symmetry inherented to
the Fock representation only. We would like to stress that this additional
relations (4)-(5) are not valid in other representations in which the
central element $\zeta $ takes the values different from zero. This means
that including of the relation (5) (or, equivalent, relations (4)) into the
list of defining commutation relations (3) drastically restricts possible
representations to the Fock representation only. It can be shown (see
[26],for example) that such extension of the commutation relations (3)
generates the algebra {\cal A$_{q,1/q}$} isomorphic to the quantum algebra U$
_qsl(2).$ Let us remark that there is one more important difference between $
q$-oscillator algebra {\cal A $_q$} and its restricted form {\cal A$_{q,1/q}$
}. Namely, the {\cal A$_{q,1/q}$} may be endowed with a Hopf algebra
structure  inherited from quantum algebra U$_qsl(2)$ one, but such structure
for {\cal A$_q$}  is still unknown and probably doesn't exist at all.

{\bf 4}. The first example of $q$-bosonization of quantum group $GL_q(2)$
was suggested in [13] (the early results [10-12] concern with the specific
particular case ''$q$ is a root of unity''). The example given in [13] i is
based on the isomorphism between Fock space for $q$-oscillator and $q$
-analogue f a Verma module (constructed in [13]) in which representation of
the quantum group $GL_q(2)$ acts. As a result the authors received two
following variants of $q$-bosonization of the quantum group $GL_q(2)$
$$
T=\left(
\begin{array}{cc}
a & b \\
c & d
\end{array}
\right) =\left(
\begin{array}{cc}
\lambda \mu \nu q^Na_{-} & \mu q^N \\
\nu q^N & a_{+}
\end{array}
\right) ,
$$
$$
T_1=\left(
\begin{array}{cc}
a & b \\
c & d
\end{array}
\right) =\left(
\begin{array}{cc}
-\lambda \mu \nu a_{+} & \mu q^{-N} \\
\nu q^{-N} & q^{-N}a_{-}
\end{array}
\right) ,
$$
where, as usual, $\lambda =q-q^{-1}$ and $\mu ,\nu $ are parameters of the
realization.

It is also possible, of course, to consider another variants of $q$
-bosonization for $GL_q(2)$. For example, realization, using only creation
operator $a_{+}$, has the form
$$
T_2=\left(
\begin{array}{cc}
a & b \\
c & d
\end{array}
\right) =\left(
\begin{array}{cc}
\mu q^N & a_{+} \\
a_{+} & \mu ^{-1}q^{-N}({\cal D}_q+qa_{+}^2)
\end{array}
\right) .
$$
Another example using two independent (mutually commuting) $q$
-oscillators $\{a_{\pm }^{(i)},N_i;i=1,2\}$ has the form
$$
T_3=\left(
\begin{array}{cc}
a & b \\
c & d
\end{array}
\right) =\left(
\begin{array}{cc}
\lambda \mu a_{+}^{(1)}a_{-}^{(1)}q^{N_2-1} & \mu \nu \sigma ^{-1}X_2 \\
\sigma q^{N_2} & \nu W_1^{-1}
\end{array}
\right) ,
$$
where X$_2=\lambda
a_{+}^{(2)}a_{-}^{(2)}+q^{-N_2},\;W_1=qa_{+}^{(1)}a_{-}^{(1)}+q^{-N_1}$ and W
$_1^{-1}$ is understand in the sense of the formal power series. In this
case ${\cal D}_q=-\mu \nu q^{-1}.$

{\bf 5}.  We'll  show now, that the Gauss decomposition of a quantum matrix $
T$ allows us to give the rather simple procedure for construction of $q$
-bosonization for $GL_q(2),$ which has more or less natural generalization
not only to the case of quantum groups $GL_q(n)$ with $n>2,$ but also for
uantum groups of the series $B_n,C_n$ and $D_n$ [21,22]. Let us return to
the simplest case of the quantum group $GL_q(2)$ and consider the Gauss
decomposition of it's quantum matrix $T=\left(
\begin{array}{cc}
a & b \\
c & d
\end{array}
\right) :$
\begin{equation}
T=T_LT_DT_R=\left(
\begin{array}{cc}
1 & u \\
0 & 1
\end{array}
\right) \left(
\begin{array}{cc}
A & 0 \\
0 & B
\end{array}
\right) \left(
\begin{array}{cc}
1 & 0 \\
z & 1
\end{array}
\right) =\left(
\begin{array}{cc}
A+uBz & uB \\
Bz & B
\end{array}
\right) .
\end{equation}
Such Gauss decomposition was firstly considered in the work [20] (remark,
that such decomposition is  used also in [23-25]) where it was noticed that
such decomposition gives the new set of generators for $GL_q(2),$ which have
more simple commutation rules than relations (1) for original generators .
Indeed, this relations have the following ($q$-Weyl) form
\begin{equation}
AB=BA,\;Au=quA,\;Az=qzA,\;uB=qBu,\;zB=qBz,\;uz=zu.
\end{equation}
If we suppose the invertibility of the initial generator $d$ in $GL_q(2)$
(or adding $d^{-1}$ to the list of generators)  we have
$$
B=d,\quad z=d^{-1}c,\quad u=bd^{-1},\quad A=a-bd^{-1}c,
$$
from which it follows that
$$
T=\left(
\begin{array}{cc}
A+uBz & uB \\
Bz & B
\end{array}
\right) =\left(
\begin{array}{cc}
1 & bd^{-1} \\
0 & 1
\end{array}
\right) \left(
\begin{array}{cc}
a-bd^{-1}c & 0 \\
0 & d
\end{array}
\right) \left(
\begin{array}{cc}
1 & 0 \\
d^{-1}c & 1
\end{array}
\right) .
$$
The element $d^{-1}$ has the following commutation relations
$$
d^{-1}c=qcd^{-1},\;d^{-1}b=qbd^{-1},\;dd^{-1}=d^{-1}d,\;
$$
$$
d^{-1}a-q^2ad^{-1}=(1-q^2){\cal D}_q(d^{-1})^2,\;{\cal D}_qd^{-1}=d^{-1}
{\cal D}_q.
$$
Note, that quantum determinant has the same value as before and is equal to
the usual determinant of central diagonal matrix of the Gauss decomposition
(6)
$$
{\cal D}_q=ad-qbc=(a-bd^{-1}c)d=AB=\det T_D.
$$
Let us remark also, that if we take another variant of the Gauss
decomposition
$$
T=\left(
\begin{array}{cc}
a & b \\
c & d
\end{array}
\right) =\left(
\begin{array}{cc}
1 & 0 \\
u & 1
\end{array}
\right) \left(
\begin{array}{cc}
A & 0 \\
0 & B
\end{array}
\right) \left(
\begin{array}{cc}
1 & z \\
0 & 1
\end{array}
\right) =\left(
\begin{array}{cc}
A & Az \\
uA & uAz+B
\end{array}
\right) ,
$$
instead of (6), then we receive the same commutation rules (7) for modified
generators as before.

{\bf 6}. Now let us turn to the problem of $q$-bosonization. Because of the
commutation rules for the new generators are simplier then fundamental
relations (1) it is natural to begin with them. We start with the case of
two independent (mutually commuting) $q$-oscillators $\{a_{\pm
}^{(i)},N_i;i=1,2\}.$ Comparing the commutation rules (3) and (7) we suppose
that
$$
\begin{array}{c}
u=\Psi _u(a_{\pm }^{(1)},N_1)a_{\pm }^{(1)};\quad z=\Psi _z(a_{\pm
}^{(2)},N_2)a_{\pm }^{(2)}; \\
\vspace{.2cm} \\ A=\Psi _A(a_{\pm }^{(i)},N_i);\quad B=\Psi _B(a_{\pm
}^{(i)},N_i);
\end{array}
$$
Here functions of the type $\Psi _u,$ etc. are understanding as formal
series or polynomials in its arguments and {\it a priori }are rather
arbitrary. Under such assumptions the commutation relations $[u,z]=0=[A,B]$
are fulfilled automatically, but other $q$-commutators give some additional
restrictions. For example, for $u\sim a_{+}^{(1)}$ commutator $Au=quA$ reads
\begin{equation}
\Psi _A(a_{\pm }^{(i)},N_i)a_{+}^{(1)}=qa_{+}^{(1)}\Psi _A(a_{\pm
}^{(i)},N_i).
\end{equation}
{}From this we received (after expanding $\Psi _A$ into the series and moving $
a_{+}^{(1)}$ to the right side)
$$
a_{+}^{(1)}\Psi _A(a_{+}^{(i)},a_{-}^{(i)},N_i)=\Psi
_A(q^{-1}a_{+}^{(1)}a_{-}^{(1)}-q^{-N_1},a_{+}^{(2)},a_{-}^{(2)},N_1-1,N_2)a_{+}^{(1)}
$$
If we insert the last expression into the former one we arrived to the
relation
\begin{equation}
\Psi _A(a_{+}^{(1)},a_{-}^{(1)},N_1)=q\Psi
_A(q^{-1}a_{+}^{(1)}a_{-}^{(1)}-q^{-N_1},N_1-1)a_{+}^{(1)}
\end{equation}
where we omit variables related to the second oscillator which are not
changing in such processes. After repeating this procedure with all other
commutators from (7) we received
\begin{eqnarray}\Psi _A(a_{+}^{(i)},a_{-}^{(i)},N_i)=q\Psi _A(q^{\mp
1}a_{+}^{(i)}a_{-}^{(i)}\mp q^{-N_i},N_1\mp 1,N_2) \\
\Psi _B(a_{+}^{(i)},a_{-}^{(i)},N_i)=q^{-1}\Psi _B(q^{\mp
1}a_{+}^{(i)}a_{-}^{(i)}\mp q^{-N_i},N_1\mp 1,N_2)
\end{eqnarray}Let us note that $q$-commutators does not limit choice of the
functions $\Psi _u,\Psi _z.$ Although we can not give general solution of
the equations (8)-(11) the concrete particular solutions can be finded
without difficulty. For example the choice
$$
\Psi _{u\;}=\alpha ,\;\Psi _z=\beta ,\;\Psi _A=\gamma q^{N_1-N_2},\;\Psi
_B=\delta q^{N_2-N_1}
$$
gives
\begin{equation}
u=\alpha a_{+}^{(1)},\;z=\beta a_{-}^{(2)},\;A=\gamma q^{N_1-N_2},\;B=\delta
q^{N_2-N_1};
\end{equation}
and defines the following realization for the matrix $T$ of $GL_q(2)$
-generators
$$
\begin{array}{c}
T=\left(
\begin{array}{cc}
a & b \\
c & d
\end{array}
\right) =\left(
\begin{array}{cc}
\gamma q^{N_1-N_2}+\alpha \beta \delta q^{N_2-N_1}a_{+}^{(1)}a_{-}^{(2)} &
\alpha \delta a_{+}^{(1)}q^{N_2-N_1} \\
\beta \delta q^{N_2-N_1}a_{-}^{(2)} & \delta q^{N_2-N_1}
\end{array}
\right) ; \\
\vspace{.2cm} \\ {\cal D}_q=\det _qT=AB=\gamma \delta .
\end{array}
$$
Slightly more general choice
$$
u=\alpha a_{+}^{(1)},\;z=\beta a_{-}^{(2)},\;A=\gamma X_1Y_2,\;B=\delta
Y_1X_2,
$$
where
$$
X_i=\lambda a_{+}^{(i)}a_{-}^{(i)}+q^{-N_i},\quad Y_i=\lambda
a_{+}^{(i)}a_{-}^{(i)}-q^{-N_i+1},\quad \lambda =q-q^{-1},
$$
gives
$$
\begin{array}{c}
T=\left(
\begin{array}{cc}
a & b \\
c & d
\end{array}
\right) =\left(
\begin{array}{cc}
\gamma X_1Y_2+\alpha \beta \delta qY_1X_2a_{+}^{(1)}a_{-}^{(2)} & \alpha
\delta a_{+}^{(1)}Y_1X_2 \\
\beta \delta X_1Y_2a_{-}^{(2)} & \delta Y_1X_2
\end{array}
\right) ; \\
\vspace{.2cm} \\ {\cal D}_q=\det _qT=AB=\gamma \delta X_1X_2Y_1Y_2.
\end{array}
$$

There is also a one $q$-boson realization. In particular, the above
mentioned variant from [13] can be obtained if we take
$$
u=\mu q^NW^{-1}a_{-},\quad z=\nu qW^{-1}q^Na_{-},
$$
$$
A=\mu \nu (\lambda -q^{N+1}W^{-1})q^Na_{-},\quad B=a_{+},
$$
where $W=qa_{+}a_{-}+q^{-N}.$

{\bf 7}. The developed above procedure of construction of $q$-bosonization
is pure algebraic one and not based on concrete choice of representation of
quantum group as well as $q$-deformed oscillator. At the same time it allows
using concrete representations of $q$-oscillator algebras to build the
quantum group representations in respected spaces. Let us consider the Fock
representation for two independent $q$-oscillators with the basis
$$
|n,m>=|n>\otimes
|m>=([n]![m]!)^{-1/2}(a_{+}^{(1)})^n(a_{+}^{(2)})^m|0>\otimes |0>
$$
in  the Fock space {\cal H}$\otimes ${\cal H}. Then for realization (12) we
have
$$
u|n,m>=\alpha \sqrt{\left[ n+1\right] }|n+1,m>,\quad z|n,m>=\beta \sqrt{
\left[ m\right] }|n,m-1>,
$$
$$
A|n,m>=\gamma q^{n-m}|n,m>,\quad B|n,m>=\delta q^{m-n}|n,m>.
$$

These relations define the following representation of the quantum group $
GL_q(2)$
$$
a|n,m>=(A+uBz)|n,m>=
$$
$$
\gamma q^{n-m}|n,m>+\alpha \beta \delta q^{m-n-1}|n,m>\sqrt{\left[
n+1\right] [m]}|n+1,m-1>,
$$
$$
b|n,m>=uB|n,m>=\beta \delta q^{m-n}\sqrt{\left[ n+1\right] }|n+1,m>,
$$
$$
c|n,m>=Bz|n,m>=\beta \delta q^{m-n-1}\sqrt{\left[ m\right] }|n,m-1>,
$$
$$
d|n,m>=B|n,m>=\delta q^{m-n}|n,m>.
$$
Let us consider the well-known realization (see f.ex. [6]) of the Fock
representation for $q$-oscillator by $q$-difference derivative in the space
{\cal P} of polynomials $p(w)$ on variable $w$ with basis $
\{w^n\}_{n=0}^\infty $ :
$$
\begin{array}{c}
p(w)a_{+}\equiv M_w\,:p(w)\longmapsto wp(w);\;w^n\longmapsto w^{n+1}, \\
\vspace{.2cm} \\ a_{-}\equiv ^qD_w\,\,:\,p(w)\longmapsto ^qD_w\,\,p(w)=
\frac{p(qw)-p(q^{-1}w)}{(q-q^{-1})w};\;w^n\longmapsto [n]w^{n-1}, \\
\vspace{.2cm} \\ q^N\equiv {}^qK_w\,\,:\,p(w)\longmapsto
{}^qK_w\,\,p(w)=p(qw);\;w^n\longmapsto q^nw^n, \\
\vspace{.2cm} \\ N\equiv w\frac d{dw}\,:\,p(w)\longmapsto wp^{\prime
}(w)=p(qw);\;w^n\longmapsto nw^n.
\end{array}
$$
Then for the same $q$-bosonization (11) of the quantum group $GL_q(2)$ we
have the following realization of generators on the space {\cal P}$\otimes $
{\cal P} by $q$-difference derivatives
$$
\begin{array}{c}
a=\gamma (^qK_w)(^qK_v)^{-1}+\alpha \beta \delta M_w(^qK_w)^{-1}(^qK_v)(^q
\negthinspace \negthinspace {}D_v), \\ \vspace{.2cm} \\ b=\alpha \delta
M_w(^qK_w)^{-1}(^qK_v),\quad c=\beta \delta (^qK_w)^{-1}(^qK_v)(^q
\negthinspace \negthinspace {}D_v),\quad d=\delta (^qK_w)^{-1}(^qK_v).
\end{array}
$$
Analogous formulae, of course, are valid and for other $q$-bosonizations of
the quantum group $GL_q(2)$.

The authors thanks Prof. P.P.Kulish for useful discussions. The work of one
of us (D.E.V.) was supported in part by the Russian Fond of Fundamental
Researches (grant No 94-01-01157-a).


\end{document}